# TORIC ARC SCHEMES AND $q$-ENUMERATION OF LATTICE POINTS

DAVID ANDERSON AND ANIKET SHAH

ABSTRACT. We introduce a natural weighted enumeration of lattice points in a polytope, and give a Brion-type formula for the corresponding generating function. The weighting has combinatorial significance, and its generating function may be viewed as a generalization of the Rogers-Szegő polynomials. It also arises from the geometry of the toric arc scheme associated to the normal fan of the polytope. We show that the asymptotic behavior of the coefficients at $q = 1$ is Gaussian.

## 1. INTRODUCTION

For each lattice polytope $P$ there is a divisor $D$ on a projective toric variety, so that the lattice points in $P$ are in natural bijection with a basis of the global sections of the line bundle $\mathcal{O}(D)$. This correspondence forms the foundation of a bridge relating algebraic geometry and combinatorics that has been exploited extensively in many contexts in the past several decades. Brion's formula is a striking example: it expresses the generating function of lattice points in $P$ (a Laurent polynomial) as a sum of rational functions determined by the tangent cones at each vertex of $P$. The simplest case is when $P$ is the interval from 0 to $m$, in which case Brion's formula is the elementary identity

$$(*) \qquad 1 + x + x^2 + \cdots + x^m = \frac{1}{1-x} + \frac{x^m}{1-x^{-1}}.$$

In this article, we establish a $q$-analogue of Brion's formula. To state the formula, we require the *$q$-Pochhammer symbol*

$$(x;q)_d = \prod_{i=1}^{d}(1 - xq^{i-1}),$$

_________
*Date*: February 7, 2023.
The authors were partially supported by NSF CAREER DMS-1945212. The second author was also supported by Charles University project PRIMUS/21/SCI/014.

as well as its infinite limit $(x;q)_\infty$, which is an analytic function when $|q| < 1$. The *q-multinomial coefficients* are polynomials in $q$ defined by

$$\begin{bmatrix} m \\ k_1, \ldots, k_r \end{bmatrix}_q = \frac{(q;q)_m}{(q;q)_{k_1} \cdots (q;q)_{k_r}},$$

where $m = k_1 + \cdots + k_r$. They specialize to 1 at $q = 0$ and to the usual multinomial coefficients at $q = 1$. When $r = 2$, these are the *q-binomial coefficients* and one usually writes $\begin{bmatrix} m \\ k \end{bmatrix}_q$ in place of $\begin{bmatrix} m \\ k, m-k \end{bmatrix}_q$.

Here is our $q$-analogue of $(*)$:

$(*_q)$

$$\begin{bmatrix} m \\ 0 \end{bmatrix}_q + \begin{bmatrix} m \\ 1 \end{bmatrix}_q x + \begin{bmatrix} m \\ 2 \end{bmatrix}_q x^2 + \cdots + \begin{bmatrix} m \\ m \end{bmatrix}_q x^m$$

$$= \frac{(q;q)_m}{(q;q)_\infty} \left( \frac{1}{(x;q)_\infty} \sum_{d \geq 0} \frac{q^{md}}{(q^{-1};q^{-1})_d (xq^{-1};q^{-1})_d} + \frac{x^m}{(x^{-1};q)_\infty} \sum_{d \geq 0} \frac{q^{md}}{(q^{-1};q^{-1})_d (x^{-1}q^{-1};q^{-1})_d} \right).$$

With some algebraic manipulation, the Jacobi triple product identity can be deduced from the case $m = 0$.

Our main theorem is a generalization of $(*_q)$ to higher-dimensional polytopes. Let $M$ be a lattice of rank $n$, so $M \cong \mathbb{Z}^n$, and let $N = \mathrm{Hom}(M,\mathbb{Z})$ be the dual lattice. The corresponding real vector spaces are $M_\mathbb{R}$ and $N_\mathbb{R}$. We write a polytope $P \subseteq M_\mathbb{R}$ as an intersection of half-spaces, so

$$P = \bigcap_{i=1}^r \{u \,|\, \langle u, v_i \rangle \geq -a_i\}$$

where $v_i \in N$ is the primitive inward normal vector defining the $i^{\text{th}}$ facet of $P$, and $a_i$ is an integer.

The inward normal vectors $v_1, \ldots, v_r$ determine a homomorphism $\mathbb{Z}^r \to N$. We denote the kernel by $A$, writing $\beta \colon A \to \mathbb{Z}^r$ for the inclusion and $A_+ = \beta^{-1}(\mathbb{Z}^r_{\geq 0})$ for the intersection with the positive orthant.

The result is cleanest under some assumptions on the polytope $P$. We assume $P$ is *smooth*: each vertex is contained in exactly $n$ facets, and the corresponding primitive normal vectors form a basis for $N$. We further assume the polytope is *radially symmetric*, meaning that $\sum_{i=1}^r v_i = 0$. These conditions can be removed, at the expense of adding extra inequalities $\langle u, v \rangle \geq -a$; see Remark 4.5. On the other hand, there are well-known families of polytopes which satisfy the stated hypotheses—for example, the (generalized) permutahedra are smooth and radially symmetric.

For each vertex $p$ of $P$, let $I(p) \subseteq \{1, \ldots, r\}$ be the set of indices $i$ so that $v_i$ is an inward normal vector of a facet containing $p$, so $\{v_i \mid i \in I(p)\}$ is a basis for $N$. Let $\{u_i(p) \mid i \in I(p)\}$ be the dual basis of $M$. The vector $u_i(p)$ is the primitive vector along the edge of $P$ which is not contained in the facet defined by $v_i$. For each $d \in A_+$ and vertex $p$ of $P$, we define

$$\mathsf{J}_{d,p} = \left( \prod_{i \in I(p)} \frac{1}{(x^{u_i(p)} q^{-1}; q^{-1})_{\beta(d)_i}} \right) \left( \prod_{j \notin I(p)} \frac{1}{(q^{-1}; q^{-1})_{\beta(d)_j}} \right).$$

**Theorem 1.** *Let $P$ be a smooth, radially symmetric polytope, realized as an intersection of $r$ half-spaces $\{u \mid \langle u, v_i \rangle \geq -a_i\}$, with notation as above. Then*

$$\sum_{u \in P \cap M} \left[ \begin{array}{c} |a| \\ \langle u, v_1 \rangle + a_1, \ldots, \langle u, v_r \rangle + a_r \end{array} \right]_q x^u = \frac{(q;q)_{|a|}}{(q;q)_\infty^{r-n}} \sum_{p \in V(P)} \sum_{d \in A_+} \frac{x^p q^{\sum a_i \beta(d)_i} \cdot \mathsf{J}_{d,p}}{\prod_{i \in I(p)} (x^{u_i(p)}; q)_\infty},$$

*where $|a| = a_1 + \cdots + a_r$ and $V(P)$ is the set of vertices of $P$.*

The left-hand side belongs to $\mathbb{Z}[q][M]$; i.e., it is a Laurent polynomial in $x$ and a polynomial in $q$. The right-hand side is a formal power series in both $x$ and $q$, and it seems (to us) surprising that a $q \to 1$ limit exists, because each term has a pole of infinite order at $q = 1$. The $q \to 0$ limit, on the other hand, recovers Brion's formula in the smooth, radially symmetric case. (To see this, note that $\mathsf{J}_{d,p}$ is 0 at $q = 0$ unless $d = 0$, in which case it is 1.)

The theorem follows from Theorem 4.4, which relaxes the requirement of radial symmetry and gives a similar formula for any smooth polytope $P$. As in [Br88], the basic proof technique is localization in the equivariant K-theory of a toric variety. Our main observation is that a canonical $q$-enumeration is provided by *toric quasimap spaces*, introduced by Morrison-Plesser and Givental in the context of Gromov-Witten theory [MP95, G98]. These spaces fit together to form an ind-variety contained in the *toric arc scheme* studied by Arkhipov-Kapranov [AK06]. This infinite-dimensional scheme provides the geometric context for the infinite products appearing in the theorem.

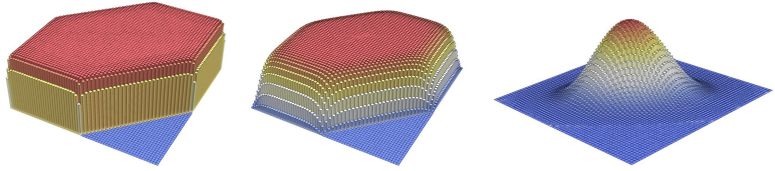

FIGURE 1. Approximations to the loop-DH measure. The polytope $P$ is a lattice hexagon, the convex hull of $(0,0)$, $(1,0)$, $(0,1)$, $(2,1)$, $(1,2)$, and $(2,2)$. The three images depict the integer-point transforms of the LHS of the formula in Theorem 1, for the dilation $kP$ with $k = 30$, evaluated at $q = 0.2$, $q = 0.6$, and $q = 0.9$. (The vertical scales differ.) Images generated by Mathematica.

A different $q$-enumeration of lattice points was investigated by Chapoton, who introduced a $q$-analogue of the Ehrhart polynomial [Ch16]. His $q$-enumeration arises by specializing each Laurent monomial $x^u$ to a power $q^{l(u)}$, for some linear function $l$ on the lattice.

The left-hand side of ($*_q$) is the $n^{\text{th}}$ *Rogers-Szegő polynomial*. This is a classical family of polynomials, which were shown to be orthogonal with respect to a certain measure on the circle by Szegő in 1926 [Sz26]. Multivariate versions of them play a role in representation theory and combinatorics [H97, V12, CV14]. In Section 5, we show how the generalized Rogers-Szegő polynomials, associated to more general polytopes, are related by $q$-difference operators.

We also consider a natural analogue of Duistermaat-Heckman (DH) measure for the quasimap ind-variety. Classical DH measure is the pushforward of Liouville measure along the moment map of a symplectic manifold with torus action, but this requires modification for infinite-dimensional spaces. Adapting the construction of DH measure from [BP90] to the ind-variety, we define a loop-Duistermaat Heckman probability measure by changing the rescaling used in a limiting procedure. Our second result, which describes $q = 1$ behavior of the coefficients of the generalized Rogers-Szegő polynomial of $kP$ as $k$ goes to infinity, is that the loop-DH measure is Gaussian. We denote by $m_D$ a particular distinguished point inside $P$ (see Definition 6.6).

**Theorem 2.** *Let $P$ be smooth and radially symmetric, and let $X$ be the associated toric variety. The polytope $P$ determines very ample divisors $D$ on $X$ and $D^0$ on the quasimap ind-variety $\mathcal{Q}(X)_\infty$. The loop-DH probability measure of $\mathcal{Q}(X)_\infty$ with respect $D^0$ is*

*Gaussian with mean $m_D$, and with covariance matrix associated to the quadratic form*

$$u \mapsto \sum_{i=1}^{r} \frac{1}{\langle m_D, v_i \rangle + a_i} \langle u, v_i \rangle^2.$$

In other words, the loop-DH measure has density function given by a multiple of

$$u \mapsto \exp\left(-\frac{1}{2} \sum_{i=1}^{r} \frac{\langle u - m_D, v_i \rangle^2}{\langle m_D, v_i \rangle + a_i}\right).$$

*Acknowledgements.* We thank Hsian-Hua Tseng for helpful conversations about quasimaps and the J-function, as well as Frédéric Chapoton, Christian Haase, and Benjamin Nill for correspondence about lattice points in polytopes.

## 2. Toric varieties

By way of fixing notation and conventions, we quickly review basic facts about toric varieties. Our main references are [CLS11] and [F93]. We work over the field $\mathbb{C}$ (which for our purposes may be replaced by any field). Let $T$ be an $n$-dimensional torus with character group $M = \text{Hom}(T, \mathbb{C}^*)$ and co-character group $N = \text{Hom}(\mathbb{C}^*, T)$, so $N = \text{Hom}(M, \mathbb{Z})$ is dual to $M$, and both $N$ and $M$ are (non-canonically) isomorphic to $\mathbb{Z}^n$. As before we write $N_\mathbb{R} = N \otimes_\mathbb{Z} \mathbb{R}$ and $M_\mathbb{R} = M \otimes_\mathbb{Z} \mathbb{R}$ for the corresponding real vector spaces.

Let $\Delta$ be a rational polyhedral fan in $N_\mathbb{R}$, i.e. a collection of rational polyhedral convex cones which fit together along their faces. We write $\Delta^{(k)} \subseteq \Delta$ for the subset of $k$-dimensional cones. Of particular importance are the rays $\Delta^{(1)} = \{\rho_1, \ldots, \rho_r\}$. Each ray $\rho_i$ has a primitive generator $v_i \in N$.

We assume $\Delta$ is smooth and complete. So there is an exact sequence

$$0 \to A \xrightarrow{\beta} \mathbb{Z}^r \to N \to 0,$$

where $\mathbb{Z}^r \to N$ sends the $i^{\text{th}}$ standard basis vector to the primitive generator of the $i^{\text{th}}$ ray, $e_i \mapsto v_i$. The kernel $A$ is isomorphic to $\mathbb{Z}^{n-r}$.

Dualizing, one has an exact sequence

$$0 \to M \to \mathbb{Z}^r \to B \to 0,$$

where $B = \text{Hom}(A, \mathbb{Z})$. Let $G \subseteq (\mathbb{C}^*)^r$ be the subtorus corresponding to the surjection of lattices $\mathbb{Z}^r \to B$, so we have $G \cong (\mathbb{C}^*)^{n-r}$ and $T = (\mathbb{C}^*)^r / G$.

We will make use of the *Cox construction*, which realizes the toric variety $X(\Delta)$ as a GIT quotient of $\mathbb{C}^r$ by $G$. Each subset $I \subseteq \{1, \ldots, r\}$ determines a coordinate subspace $E_I = \{e_i^* = 0 \,|\, i \in I\} \subseteq \mathbb{C}^r$, as well as a collection of rays. Let $Z(\Delta) \subseteq \mathbb{C}^r$ be the union of those coordinate subspaces $E_I$ such that the corresponding set of rays $\{\rho_i \,|\, i \in I\}$ is not contained in any cone of $\Delta$. Then

$$X(\Delta) = (\mathbb{C}^r \smallsetminus Z(\Delta))/G.$$

Basic facts from toric geometry say that

$$B \cong \mathrm{Pic}(X) \cong H^2(X, \mathbb{Z})$$

and the cone of effective divisors is the image of $\mathbb{Z}_{\geq 0}^r$ in $B$. Dually,

$$A \cong H_2(X, \mathbb{Z}),$$

and the cone of nef curves is the preimage of $\mathbb{Z}_{\geq 0}^r$ under the embedding $\beta \colon A \hookrightarrow \mathbb{Z}^r$. This cone is written $A_+ \subseteq A$.

When $X(\Delta)$ is projective, an ample line bundle $\mathcal{O}(D)$ corresponds to a polytope $P$ with normal fan $\Delta$. Facets of $P$ correspond to $T$-invariant divisors of $X$, and vertices of $P$ correspond to fixed points. We sometimes abuse notation by identifying vertices and fixed points, writing both $p \in P$ and $p \in X^T$.

Sections of $\mathcal{O}(D)$ may be described in terms of the Cox construction as follows. Suppose $P$ is defined as the intersection of half-spaces $\{u \,|\, \langle u, v_i \rangle \geq -a_i\}$, for $i = 1, \ldots, r$; then $D = \sum a_i D_i$, where $D_1, \ldots, D_r$ are the $T$-invariant divisors. Let $f_1, \ldots, f_r$ be standard coordinates on $\mathbb{C}^r$, generating the *Cox ring* of $X$. Then, for nonnegative integers $b_i$, a Laurent monomial

$$x^u = \prod_{i=1}^{r} f_i^{b_i - a_i}$$

is a $T$-equivariant section of $\mathcal{O}(D)$ if and only if $\sum (b_i - a_i) D_i = 0$ in $\mathrm{Pic}(X)$; equivalently, $u = \sum (b_i - a_i) e_i$ lies in $M \subset \mathbb{Z}^r$, where $e_1, \ldots, e_r$ is the standard basis. By construction, $\langle u, v_i \rangle = b_i - a_i \geq a_i$ for each $i$, so $u \in P$.

## 3. Quasimap spaces and arc schemes

Next we review the construction of the toric quasimap space; see [MP95, G98] for proofs and details. We are interested in parametrizing maps $f \colon \mathbb{P}^1 \to X(\Delta)$ of *degree $d$*, that is, $f_*[\mathbb{P}^1] = d$ in $A = H_2(X)$. Let $\mathrm{Hom}_d(\mathbb{P}^1, X)$ be the space of such maps. To describe this space, one can lift such maps to $\mathbb{C}^r$, where they can be specified as an $r$-tuple of univariate polynomials. We consider only degrees $d$ lying in $A_+$. In general, $A_+$ is properly contained in the semigroup of all effective curves.

For any $r$-tuple of nonnegative integers $\delta = (\delta_1, \ldots, \delta_r)$, let

$$\mathbb{C}^r_\delta = \mathbb{C}[t]_{\leq \delta_1} \oplus \cdots \oplus \mathbb{C}[t]_{\leq \delta_r}$$

where $\mathbb{C}[t]_{\leq k} = \{f(t) = f^{(0)} + f^{(1)}t + \cdots + f^{(k)}t^k\}$ is the $(k+1)$-dimensional space of polynomials of degree at most $k$.

The vector space $\mathbb{C}^r_\delta$ has dimension $r + \delta_1 + \cdots + \delta_r$, so the torus $(\mathbb{C}^*)^{r+\sum \delta_i}$ acts coordinate-wise. The torus $(\mathbb{C}^*)^r$ embeds "diagonally", so that

$$(z_1, \ldots, z_r) \cdot (f_1(t), \ldots, f_r(t)) = (z_1 f_1(t), \ldots, z_r f_r(t)).$$

So the subtorus $G \subseteq (\mathbb{C}^*)^r$ also embeds in $(\mathbb{C}^*)^{r+\sum \delta_i}$. The quasimap space is constructed as a certain GIT quotient of $\mathbb{C}^r_\delta$ by this action of $G$.

For any subset $I \subseteq \{1, \ldots, r\}$, let $E_{I,\delta} \subseteq \mathbb{C}^r_\delta$ be the locus where $f_i \equiv 0$ for each $i \in I$. Viewing $E_I = \mathbb{C}^{r-\#I}$ as a coordinate subspace and writing $\delta(\widehat{I})$ for the subsequence of $\delta$ omitting $i \in I$, this is the same as $\mathbb{C}^{r-\#I}_{\delta(\widehat{I})}$.

Let

$$Z(\Delta)_\delta = \bigcup_I E_{I,\delta},$$

the union over $I$ such that $\{\rho_i \mid i \in I\}$ is not contained in a cone of $\Delta$.

The space of *(toric) quasimaps* is defined as

$$\mathcal{Q}(\Delta)_d = (\mathbb{C}^r_{\beta(d)} \setminus Z(\Delta)_{\beta(d)})/G.$$

By construction, it contains $\mathrm{Hom}_d(\mathbb{P}^1, X)$ as a dense open subset.

**Proposition 3.1.** *The toric quasimap space $\mathcal{Q}(\Delta)_d$ is a toric variety of dimension $n + \sum \beta(d)_i$, smooth and complete whenever $X(\Delta)$ is.*

Let us write $T_d = (\mathbb{C}^*)^{r+\sum \beta(d)_i}/G$ for the torus acting on $\mathcal{Q}(\Delta)_d$. We are also interested in actions by various subtori. Using $(\mathbb{C}^*)^r \hookrightarrow (\mathbb{C}^*)^{r+\sum \beta(d)_i}$ as above, we have an inclusion of $T = (\mathbb{C}^*)^r/G$ in $T_d$, so the same torus acting on $X(\Delta)$ also acts on $\mathcal{Q}(\Delta)_d$. On the other hand, there is a loop rotation action of $\mathbb{C}^*$ on $\mathbb{C}^r_{\beta(d)}$ by

$$\zeta \cdot (f_1(t), \ldots, f_r(t)) = (f_1(\zeta^{-1}t), \ldots, f_r(\zeta^{-1}t)),$$

and this descends to an action of $\mathbb{C}^*$ on $\mathcal{Q}(\Delta)_d$. So there is a subtorus $\mathbb{T} = T \times \mathbb{C}^* \subseteq T_d$ acting on $\mathcal{Q}(\Delta)_d$. Localization with respect to $\mathbb{T}$ will be the primary tool in proving the main theorem.

**Proposition 3.2.** *The fixed locus $(\mathcal{Q}(\Delta)_d)^{\mathbb{C}^*}$ by the loop rotation action is a union of components $\mathcal{Q}(\Delta)_d^{(d')}$, for each decomposition $d' \leq d$ (meaning $d = d' + d''$ is a decomposition as a sum of nef classes). Furthermore, there is a T-equivariant isomorphism $\mathcal{Q}(\Delta)_d^{(d')} \cong X(\Delta)$ for all $d'$.*

The quasimap spaces $\mathcal{Q}(\Delta)_d$ form a system of closed embeddings for $d$ varying over $A_+ \subset A = H_2(X, \mathbb{Z})$. Namely, given two such curve classes $d' \leq d$, there is a closed embedding $\mathcal{Q}(\Delta)_{d'} \hookrightarrow \mathcal{Q}(\Delta)_d$ induced by the inclusion $\mathbb{C}^r_{\beta(d')} \hookrightarrow \mathbb{C}^r_{\beta(d)}$. We refer to the limiting ind-variety as the toric polynomial space $\mathcal{Q}(\Delta)_\infty$.

The ind-variety $\mathcal{Q}(\Delta)_\infty$ sits inside an even larger space, the toric arc scheme of Arkhipov-Kapranov [AK06]. The toric arc scheme is constructed in nearly the same manner as $\mathcal{Q}(\Delta)_d$. Let

$$\mathbb{C}^r_\infty = \mathbb{C}[[t]] \oplus \cdots \oplus \mathbb{C}[[t]].$$

This space has an infinite-dimensional torus action given by scaling coordinates of the tuple of power series. For $I \subset \{1, \ldots, r\}$, let the locus $E_{I,\infty} \subseteq \mathbb{C}^r_\infty$ be defined as the tuples where all the coefficients of $f_i$ vanish for all $i \in I$, and

$$Z(\Delta)_\infty = \bigcup_I E_{I,\infty},$$

where once again the union is over $I$ such that $\{\rho_i | i \in I\}$ is not contained in a cone of $\Delta$. The Arkhipov-Kapranov toric arc scheme is

$$\Lambda^0 X = (\mathbb{C}^r_\infty \setminus Z(\Delta)_\infty)/G.$$

Just as with quasimap spaces, the torus $\mathbb{T}$ acts on $\Lambda^0 X$. The quasimap space $\mathcal{Q}(\Delta)_d$ embeds as a finite-dimensional invariant subvariety in $\Lambda^0 X$, and the following diagram commutes:

$$\begin{array}{ccc} & \mathcal{Q}(\Delta)_\infty & \\ \nearrow & & \searrow \\ \mathcal{Q}(\Delta)_d & \hookrightarrow & \Lambda^0(X). \end{array}$$

The closed points of the ind-variety $\mathcal{Q}(\Delta)_\infty$ (or equivalently, the union of closed points of $\mathcal{Q}(\Delta)_d$ over all $d \in A_+$) correspond to the tuples of power series inside $\Lambda^0 X$ where only finitely many coefficients of any given power series are non-zero. In other words, the toric polynomial space $\mathscr{Q}(\Delta)_\infty$ embeds into $\Lambda^0 X$ as the subset consisting of power series which are in fact polynomial.

**Definition 3.3.** Let $D_i$ be the $T$-invariant divisor corresponding to the ray $\rho_i$ in $X$. For each $k$ such that $\beta_i(d) \geq k \geq 0$ there are divisors $D_i^k$ in $\mathscr{Q}(\Delta)_d$, $\Lambda^0 X$, and $\mathscr{Q}(\Delta)_\infty$ defined by the vanishing of the $k^{\text{th}}$ coefficient of $f_i(t)$. For each $d \in A_+$ and each $T$-invariant divisor $D = \sum_i a_i D_i$ on $X$, we let $D^d$ be the divisor $\sum_i a_i D_i^{\beta_i(d)}$ on $\mathscr{Q}(\Delta)_\infty$. We use the same notation for analogously defined divisors on $\mathscr{Q}(\Delta)_d$ and on $\Lambda^0 X$, the base space being clear from context.

Arkhipov and Kapranov observed that $\Lambda^0 X$ admits a family of self-embeddings. Recall $\beta$ is the inclusion $H_2(X, \mathbb{Z}) \hookrightarrow \mathbb{Z}^r$ defined in Section 2. An element $d$ in the semigroup $A_+$ corresponds to a one-parameter subgroup of $G$, and by composing with the inclusion $G \hookrightarrow (\mathbb{C}^*)^r$, we can write the image of $d$ in $\mathbb{Z}^r$ explicitly as the cocharacter $(t^{\beta(d)_1}, \ldots, t^{\beta(d)_r})$ of $(\mathbb{C}^*)^r$.

**Definition 3.4.** For $d \in A_+$, let $\epsilon_d : \Lambda^0 X \to \Lambda^0 X$ be the self-embedding

$$(f_1(t), \ldots, f_r(t)) \mapsto (t^{\beta(d)_1} f_1(t), \ldots, t^{\beta(d)_r} f_r(t)).$$

This restricts to a self-embedding on the polynomial space $\mathscr{Q}(\Delta)_\infty$ which we also denote by $\epsilon_d$.

These self-embeddings commute: in fact, $\epsilon_d \circ \epsilon_{d'} = \epsilon_{d+d'}$. They are evidently equivariant with respect to the $\mathbb{T}$-action.

**Lemma 3.5.** *Let $D = \sum_i a_i D_i$ be nef on $X$. Then $D^0 = \sum_i a_i D_i^0$ restricts to a nef divisor on each $\mathscr{Q}(\Delta)_d$.*

*Proof.* We will show that $D^0$ pairs positively with any $T \times \mathbb{C}^*$-invariant effective irreducible curve $C$ in $\mathscr{Q}(\Delta)_d$.

As a toric variety, $\mathscr{Q}(\Delta)_d$ has Cox ring variables $f_i^{(j)}$, for $1 \leq i \leq r$ and $0 \leq j \leq \beta(d)_i$, and for any $k, k'$, $\frac{f_i^{(k)}}{f_i^{(k')}}$ defines a rational function on $X$ with divisor $D_i^k - D_i^{k'}$. Thus, $D^0$ is equivalent to $D^{d'} = \sum_i a_i D_i^{\beta(d')_i}$ for any $d' \leq d$ in $A_+$.

Let $\epsilon_{d'}(p)$ be a $T \times \mathbb{C}^*$-fixed point of $C$, for some $d' \leq d$ in $A_+$. If $E$ is an ample divisor on $X$, then for $s > 0$, $D + sE$ can be moved to an effective $\mathbb{Q}$-divisor $\sum_i b_i D_i$ whose support does not contain $p$. Similarly, the divisor $D^0 + sE^0$ on $\mathscr{Q}(\Delta)_d$ can be moved to $\sum_i b_i D_i^{\beta(d')_i}$, whose support does not contain $\epsilon_{d'}(p)$—consequently, its support does not contain $C$. So $C \cap (D^0 + sE^0) \geq 0$ for all $s > 0$, and therefore $C \cap D^0 \geq 0$ as well. $\square$

## 4. Lattice points and localization on the arc scheme

Now, we study the geometry of the $\mathbb{T}$-action on $\mathscr{Q}(\Delta)_\infty$, leading to the proof of our analogue of Brion's identity for $q$-weighted lattice point enumeration. To define the weighting, let $\mathscr{A} = \{(H_i, v_i, a_i)\}$ be an oriented hyperplane arrangement in $M_\mathbb{R}$. Explicitly, this means that for each $i$ we specify a primitive vector $v_i$ and integer $a_i$ such that the hyperplane $H_i$ is defined by

(1) $$\langle u, v_i \rangle = -a_i.$$

These choices determine positive and negative half-spaces $H_{i,+}$ and $H_{i,-}$, by replacing the "=" in (1) with "≥" or "≤."

Given such an arrangement $\mathscr{A}$ and any $u \in M$, we let

(2) $$g_{\mathscr{A}, u} = \prod_{i=1}^r \frac{1}{(q;q)_{\langle u, v_i \rangle + a_i}} = \prod_{i=1}^r \prod_{k=1}^{\langle u, v_i \rangle + a_i} \frac{1}{1 - q^k}.$$

Now, let $X = X(\Delta)$ be a smooth and projective toric variety, so $\Delta$ is the inward normal fan of a polytope, and let $v_1, \ldots, v_r$ once again refer to the primitive vectors in the rays of $\Delta$. Let $a_i$ be a tuple of integers further satisfying that the intersection of positive half-spaces $P = \bigcap_i H_{i,+} = \bigcap_i \{u \mid \langle u, v_i \rangle \geq -a_i\}$ is non-empty, and each hyperplane $H_i$ touches $P$. We write $\mathscr{A}_P$ for the corresponding oriented hyperplane arrangement, noting that this depends not only on $P$ but also on the rays of $\Delta$. This data determines a nef line bundle $\mathscr{O}(D) = \mathscr{O}(\sum_i a_i D_i)$ on $X$ whose global sections are in natural bijection with the lattice points in $P$. The same data also determines a line bundle $\mathscr{O}(D^0) = \mathscr{O}(\sum_i a_i D_i^0)$ on each $\mathscr{Q}(\Delta)_d$, and on $\mathscr{Q}(\Delta)_\infty$.

**Theorem 4.1.** *The equivariant Euler characteristic is given by*

(3) $$\chi_\mathbb{T}\left(\mathscr{Q}(\Delta)_\infty, \mathscr{O}(D^0)\right) = \sum_{u \in P \cap M} g_{\mathscr{A}_P, u} x^u.$$

*Proof.* Since $\mathscr{O}(D)$ is nef on $X$, the line bundle $\mathscr{O}(D^0)|_{\mathscr{Q}(\Delta)_d}$ is nef on $\mathscr{Q}(\Delta)_d$ by Lemma 3.5, so Demazure vanishing implies that

$$\chi_\mathbb{T}\left(\mathscr{Q}(\Delta)_d, \mathscr{O}(D^0)\right) = H^0\left(\mathscr{Q}(\Delta)_d, \mathscr{O}(D^0)\right)$$

as $\mathbb{T}$-modules, for each $d$. There are inclusions of $\mathbb{T}$-modules $H^0\left(\mathcal{Q}(\Delta)_d, \mathcal{O}(D^0)\right) \hookrightarrow H^0\left(\mathcal{Q}(\Delta)_{d'}, \mathcal{O}(D^0)\right)$ for $d \leq d'$, corresponding to the distinguished embeddings $\mathcal{Q}(\Delta)_d \hookrightarrow \mathcal{Q}(\Delta)_{d'}$. The infinite-dimensional $\mathbb{T}$-module $H^0\left(\mathcal{Q}(\Delta)_\infty, \mathcal{O}(D^0)\right)$, by definition, is the union of these, and $\chi_\mathbb{T}\left(\mathcal{Q}(\Delta)_\infty, \mathcal{O}(D^0)\right)$ is its graded character. To compute it, we determine the characters of each $H^0\left(\mathcal{Q}(\Delta)_d, \mathcal{O}(D^0)\right)$ and take the limit as $d \to \infty$.

As a toric variety, $\mathcal{Q}(\Delta)_d$ has Cox ring variables $f_i^{(j)}$, for $1 \leq i \leq r$ and $0 \leq j \leq \beta(d)_i$. A basis for the sections of $\mathcal{O}(D^0)$ consists of monomials $\left(\prod_i (f_i^{(0)})^{-a_i}\right) \cdot \left(\prod_i \prod_{j=0}^{\beta(d)_i} (f_i^{(j)})^{b_{i,j}}\right)$ such that the $b_{i,j}$ are non-negative, and $\sum_i \left((\sum_{j=0}^{\beta(d)_i} b_{i,j}) - a_i\right) D_i = 0$ in $\text{Pic}(X)$.

We must compute the $\mathbb{T}$-character of the span of these sections. On $X$, an element of the distinguished basis of sections of $\mathcal{O}(D)$ corresponds to a $T$-character $x^u$ where $u$ is a lattice point in $P$, or equivalently a choice of $b_i \geq 0$ such that $\sum_i (b_i - a_i) e_i \in \mathbb{Z}^r$ is the image of some $u \in M$. These choices are related via the identity $b_i = \langle u, v_i \rangle + a_i$.

With these $b_i$ fixed, consider the sections $\left(\prod_i (f_i^{(0)})^{-a_i}\right) \cdot \left(\prod_i \prod_{j=0}^{\beta(d)_i} (f_i^{(j)})^{b_{i,j}}\right)$ such that $\sum_{j=0}^{\beta(d)_i} b_{i,j} = b_i$. The character of $\prod_i f_i^{b_i - a_i}$ is $x^u$, as is that of $\prod_i (f_i^{(0)})^{b_i - a_i}$. So $\left(\prod_i (f_i^{(0)})^{-a_i}\right) \cdot \left(\prod_i \prod_{j=0}^{\beta(d)_i} (f_i^{(j)})^{b_{i,j}}\right)$ has character $x^u \prod_i \prod_{j=0}^{\beta(d)_i} (q^j)^{b_{i,j}}$. So the coefficient of $x^u$ in the graded character of $H^0\left(\mathcal{Q}(\Delta)_d, \mathcal{O}(D^0)\right)$ is

$$\sum \prod_i \prod_{j=0}^{\beta(d)_i} (q^j)^{b_{i,j}},$$

the sum over $b_{i,j} \geq 0$ such that $\sum_{j=0}^{\beta(d)_i} b_{i,j} = b_i$.

In the limit as $d \to \infty$, the upper bound disappears for the indices $j$ of $b_{i,j}$. The resulting sum is over all choices of weakly increasing sequences $0 \leq c_{i,1} \leq c_{i,2} \leq \cdots \leq c_{i,b_i}$ for each $i$ from 1 to $r$, where the summand is the statistic

$$\prod_{i=1}^r \prod_{j=1}^{b_i} q^{c_{i,j}}.$$

Holding all but one weakly increasing sequence fixed, we see that the whole sum must factor into a product over $i$:

$$\prod_{i=1}^r \sum_{0 \leq c_{i,1} \leq c_{i,2} \leq \ldots \leq c_{i,b_i}} \prod_{j=1}^{b_i} q^{c_{i,j}}.$$

But
$$\sum_{0\leq c_{i,1}\leq c_{i,2}\leq\ldots\leq c_{i,b_i}}\prod_{j=1}^{b_i} q^{c_{i,j}} = \frac{1}{(q;q)_{b_i}} = \frac{1}{(q;q)_{\langle u,v_i\rangle+a_i}}.$$

Thus the coefficient of $x^u$ in the character of $H^0(\mathcal{Q}(\Delta)_\infty, \mathcal{O}(D^0))$ is
$$\prod_{i=1}^r \frac{1}{(q;q)_{\langle u,v_i\rangle+a_i}} = g_{\mathcal{A}_P,u}.$$

This proves the theorem. □

Next we describe the $\mathbb{T}$-fixed points of $\mathcal{Q}(\Delta)_\infty$, and decompose the corresponding tangent spaces into characters of $\mathbb{T}$.

**Proposition 4.2.** *The $\mathbb{T}$-fixed points of $\mathcal{Q}(\Delta)_\infty$ are in bijection with pairs $(p,d)$, with $p \in X^T$ and $d \in A_+$.*

*Proof.* The locus fixed by loop rotation is easy to determine. By definition $\mathcal{Q}(\Delta)_\infty$ is a union of $\mathcal{Q}(\Delta)_d$ as $d$ varies over $A_+$, and the subvariety $\mathcal{Q}(\Delta)_0 \hookrightarrow \mathcal{Q}(\Delta)_\infty$ is certainly $\mathbb{C}^*$-fixed. To exhaust the $\mathbb{C}^*$-fixed components in $\mathcal{Q}(\Delta)_d$ for higher $d$, it is enough to take the disjoint union of $\epsilon_d(\mathcal{Q}(\Delta)_0)$ over $d \in A_+$, because $\epsilon_d(\mathcal{Q}(\Delta)_0) = \mathcal{Q}(\Delta)_d^{(d)} \subset \mathcal{Q}(\Delta)_d \subset \mathcal{Q}(\Delta)_\infty$. The subvariety $\mathcal{Q}(\Delta)_0$ is $T$-equivariantly isomorphic to $X$, so its $T$-fixed points are in natural bijection with those of $X$. The same holds for each $\epsilon_d(\mathcal{Q}(\Delta)_0)$, so the fixed points are simply $\epsilon_d(p)$ for $d \in A_+$ and $p \in X^T = \mathcal{Q}(\Delta)_0$. □

Recall from the introduction that for each vertex $p$ of a smooth polytope $P$, there is a subset $I(p) \subseteq \{1,\ldots,r\}$ so that the $v_i$ for $i \in I(p)$ are the inward normal vectors of the facets containing $p$, and $\{u_i(p) \mid i \in I(p)\}$ is the basis of $M$ dual to $\{v_i \mid i \in I(p)\}$. We are using the notation

$$(4) \qquad \mathsf{J}_{d,p} = \left(\prod_{i\in I(p)} \frac{1}{(x^{u_i(p)}q^{-1};q^{-1})_{\beta(d)_i}}\right)\left(\prod_{j\notin I(p)} \frac{1}{(q^{-1};q^{-1})_{\beta(d)_j}}\right).$$

**Proposition 4.3.** *Let $p$ be a fixed point in $X$, corresponding to a vertex of $P$. The equivariant multiplicity to $\mathcal{Q}(\Delta)_\infty$ at $\epsilon_d(p)$ is equal to*
$$\left(\frac{1}{(q;q)_\infty}\right)^{r-n} \frac{\mathsf{J}_{d,p}}{\prod_{i\in I(p)}(x^{u_i(p)};q)_\infty}.$$

In our usage of the term, "equivariant multiplicity" is a synonym for "Bott denominator": it is the product of factors $(1 - \chi)^{-1}$ over characters $\chi$ appearing in the cotangent space to a variety at a fixed point. The characters appearing in the denominator of $\mathsf{J}_{d,p}$ are the conormal weights to $\epsilon_d(\Lambda^0 X) \subseteq \Lambda^0 X$ at the fixed point $\epsilon_d(p)$. The other characters in the formula are cotangent weights at $p \in X = \mathcal{Q}(\Delta)_0 \subset \mathcal{Q}(\Delta)_\infty$, as one sees from the Cox description of $\mathcal{Q}(\Delta)_d$ for large enough $d$.

Now we can state and prove our main theorem:

**Theorem 4.4.** *Let P be a smooth polytope, with notation as above. Then*

$$\tag{5} \sum_{u \in P \cap M} g_{\mathcal{A}_P, u} x^u = \frac{1}{(q;q)_\infty^{r-n}} \sum_{p \in V(P)} x^p \sum_{d \in A_+} \frac{q^{\sum a_i \beta(d)_i} \cdot \mathsf{J}_{d,p}}{\prod_{i \in I(p)} (x^{u_i(p)}; q)_\infty},$$

*where $g_{\mathcal{A}_P, u}$ is the rational function defined in* (2)*, and $V(P)$ is the set of vertices of $P$.*

In the radially symmetric case, so $\sum v_i = 0$, the theorem from the introduction follows immediately by multiplying both sides by $(q;q)_{|a|}$, where $|a| = a_1 + \cdots + a_r$ as before.

*Proof.* The polytope $P$ determines a toric variety $X$, normal fan $\Delta$ and a line bundle $\mathcal{O}(\sum_i a_i D_i)$. The left-hand side is simply $\chi_{\mathbb{T}}(\mathcal{Q}(\Delta)_\infty, \mathcal{O}(\sum_i a_i D_i^0))$, by Theorem 4.1.

The right-hand side comes by applying the Atiyah-Bott formula to compute the Euler characteristic of $\mathcal{O}(\sum_i a_i D_i^0)$. This is a sum over the $\mathbb{T}$-fixed points $\epsilon_d(p)$ of the product of the $\mathbb{T}$-character of the line $i^*_{\epsilon_d(p)} \mathcal{O}(\sum_i a_i D_i^0)$ and the equivariant multiplicity at $\epsilon_d(p)$. The latter is computed by Proposition 4.3. The former is equal to $x^p \cdot q^{\sum a_i \beta(d)_i}$, as can be seen from the Cox construction of $\mathcal{Q}(\Delta)_{d'}$ for $d' \gg d$. □

**Remark 4.5.** Say a fan $\Delta$ is *radially symmetric* if the primitive generators of its rays sum to 0. The normal fan of any polytope $P$ can be refined to a smooth, radially symmetric fan $\Delta$, as follows. Given any polytope $P$, let $\Delta_0$ be its normal fan. Refine $\Delta_0$ to obtain a fan $\Delta_1$ such that for every cone $\sigma \in \Delta_1$, $-\sigma$ is also in $\Delta_1$; in particular its primitive ray generators sum to 0. Further refine $\Delta_1$ to obtain a smooth fan $\Delta$ with the same property. Now one can write $P = \bigcap \{u \langle u, v_i \rangle \geq -a_i\}$ with some redundant inequalities.

This lets us remove the radial symmetry and smoothness hypotheses on $P$. Theorem 4.4 applies to the fan $\Delta$, as does the formula of Theorem 1. However, only for smooth, radially symmetric $P$ are those formulas canonical.

**Remark 4.6.** When $X$ is a product of projective spaces, the rational function $\mathsf{J}_{d,p}$ is equal to the restriction at $p$ of the $d^{\text{th}}$ term of the J-function for quantum K-theory. More generally, it is a localization of the K-theoretic I-function; see [GT14, G15].

## 5. Generalized Rogers-Szegő polynomials and Jackson partial derivatives

In this section, we study the action of difference operators on $q$-series appearing on the left-hand side of Theorem 4.4.

Let $\Delta$ be a smooth fan, with primitive ray generators $v_1, \ldots, v_r$ summing to 0. Let $D = \sum_i a_i D_i$ be an effective $T$-divisor on the associated toric variety $X$. The condition that $D$ is effective is equivalent to the statement that the *Newton polytope* $P = \bigcap_i \{u \mid \langle u, v_i \rangle \geq -a_i\}$ contains an element of $M$. Let $|a| = \sum_i a_i$. We define the *Rogers-Szegő polynomial of D* by

$$(6) \qquad RS_D(x; q) = \sum_{u \in P \cap M} \begin{bmatrix} |a| \\ \langle u, v_1 \rangle + a_1, \ldots, \langle u, v_r \rangle + a_r \end{bmatrix}_q x^u.$$

The following comes from Theorem 4.1.

**Proposition 5.1.** *Let $D = \sum_i a_i D_i$ be a T-invariant nef divisor on a smooth complete toric variety with fan $\Delta$ and rays $v_1, \ldots, v_r$ summing to 0. Then*

$$(q; q)_{|a|} \cdot \chi_{\mathbb{T}}\left(\mathcal{Q}(\Delta)_\infty, \mathcal{O}(D^0)\right) = RS_D(x; q).$$

For generic values of $c$, the polytope $P$ corresponding to the toric divisor $D$ is exactly the Newton polytope of $RS_D(x; c)$. If the $a_i$ satisfy the further condition that each hyperplane $H_i = \{u \mid \langle u, v_i \rangle = -a_i\}$ touches $P$, we write $RS_P(x; q)$ for $RS_D(x; q)$. (Given $P$, the $a_i$ satisfying this further condition are unique.) Though it is not reflected in the notation, the polynomial $RS_P(x; q)$ depends on both $P$ and the rays of $\Delta$.

Next we prove some facts about the action of difference operators on $RS_D(x; q)$. For explicitness, we henceforth identify $N$ and $M$ with $\mathbb{Z}^n$. The $i^{\text{th}}$ $q$-shift operator $T_{i,q}$ acts on $f(x_1, \ldots, x_n)$ by

$$T_{i,q} : f(x_1, \ldots, x_i, \ldots, x_n) \mapsto f(x_1, \ldots, q x_i, \ldots, x_n),$$

and the $i^{\text{th}}$ partial $q$-derivative (or *Jackson derivative*) acts on $f(x_1, \ldots, x_n)$ by

$$\left(\frac{d}{dx_i}\right)_q : f \mapsto \frac{f - T_{i,q} f}{(1 - q) x_i}.$$

The behavior of these operators on Rogers-Szegő polynomials is best with some conditions on the polytopes and fans involved:

**Definition 5.2.** Given a fan $\Delta$ in $\mathbb{R}^n$, fix the standard basis of $\mathbb{R}^n$ and order the rays of $\Delta$ so that $v_1 = e_1, \ldots, v_n = e_n$. We say an effective divisor $D = \sum_i a_i D_i$ on $X$ is *first-orthant* if its corresponding polytope $P$ is contained in $(\mathbb{R}_{\geq 0})^n$, with $a_1 = \cdots = a_n = 0$.

For each $1 \leq i \leq n$, we denote by $P_i$ the (possibly empty) polytope in $(\mathbb{R}_{\geq 0})^n$ obtained by first replacing the condition $\{u \mid \langle u, v_i \rangle \geq 0\}$ with $\{u \mid \langle u, v_i \rangle \geq 1\}$, and then translating the result by $-e_i$. So $P_i$ is defined by $\left(\bigcap_{j=1}^{n}\{u \mid \langle u, v_j \rangle \geq 0\}\right) \cap \left(\bigcap_{j=n+1}^{r}\{u \mid \langle u + e_i, v_j \rangle \geq -a_j\}\right)$.

**Proposition 5.3.** *The partial $q$-derivatives of $RS_D(x; q)$ for $D$ a first-orthant divisor are given by*

$$\left(\frac{d}{dx_i}\right)_q RS_D(x; q) = [|a|]_q\, RS_{D + \sum_{j=n+1}^{r} \langle e_i, v_j \rangle D_j}(x; q), \tag{7}$$

*where $[b]_q = \frac{1-q^b}{1-q} = 1 + \cdots + q^{b-1}$ for an integer $b$.*

*Proof.* This is a straightforward computation:

$$\left(\frac{d}{dx_i}\right)_q RS_D(x; q) = \sum_{u \in P \cap M} \begin{bmatrix} |a| \\ \langle u, v_1 \rangle + a_1, \ldots, \langle u, v_r \rangle + a_r \end{bmatrix}_q \frac{1 - T_{i,q}}{(1-q) x_i} x^u$$

$$= \sum_{\substack{u \in P \cap M, \\ \langle u, v_i \rangle \geq 1}} \begin{bmatrix} |a| \\ \langle u, v_1 \rangle + a_1, \ldots, \langle u, v_r \rangle + a_r \end{bmatrix}_q \frac{(1 - q^{\langle u, v_i \rangle})}{(1-q)} \frac{x^u}{x_i}.$$

Replacing $u$ with $u + e_i$, this becomes

$$\sum_{u \in P_i \cap M} \begin{bmatrix} |a| \\ \ldots, \langle u, v_i + e_i \rangle, \ldots \end{bmatrix}_q \frac{(1 - q^{\langle u+e_i, v_i \rangle})}{(1-q)} x^u,$$

$$= \sum_{u \in P_i \cap M} \begin{bmatrix} |a| \\ \ldots, \langle u, v_i \rangle + 1, \ldots \end{bmatrix}_q \frac{(1 - q^{\langle u, v_i \rangle + 1})}{(1-q)} x^u,$$

noting that $a_i = 0$ by assumption. We have

$$\begin{bmatrix} |a| \\ \ldots, \langle u, v_i \rangle + 1, \ldots \end{bmatrix}_q (1 - q^{\langle u, v_i \rangle + 1})$$

$$= \begin{bmatrix} |a| - 1 \\ \ldots, \langle u, v_i \rangle, \ldots \end{bmatrix}_q (1 - q^{|a|}).$$

Thus,

$$\left(\frac{d}{dx_i}\right)_q RS_D(x;q) = \sum_{u \in P_i \cap M} \begin{bmatrix} |a| - 1 \\ \langle u, v_1 \rangle, \ldots, \langle u, v_n \rangle, \ldots, \langle u + e_i, v_r \rangle + a_r \end{bmatrix}_q \frac{(1 - q^{|a|})}{(1 - q)} x^u,$$

which is exactly $[|a|]_q RS_{D + \sum_{j=n+1}^{r} \langle e_i, v_j \rangle D_j}(x;q)$. □

**Theorem 5.4.** *Let $X$ be a smooth complete toric variety. The $\mathbb{Q}(q)$-linear span of $RS_D(x;q)$ for $D$ a first-orthant divisor on $X$ is an indecomposable representation of the (commutative) algebra spanned by $\left(\frac{d}{dx_1}\right)_q, \ldots, \left(\frac{d}{dx_n}\right)_q$.*

*Proof.* Let $D$ be any first-orthant divisor, and $P$ the Newton polytope of $D$. Let $(i_1, \ldots, i_n) \in P \cap M$ maximize the value of $x_1 + \cdots + x_n$ on $P$. Then we can calculate that $\left(\frac{d}{dx_1}\right)_q^{i_1} \cdots \left(\frac{d}{dx_n}\right)_q^{i_n} RS_D(x;q)$ is a non-zero element of $\mathbb{C}(q)$: For any other point $(j_1, \ldots, j_n)$ in $P \cap M$, there is some index $k$ such that $i_k > j_k$, so $\left(\frac{d}{dx_k}\right)_q^{i_k} x_1^{j_1} \cdots x_n^{j_n} = 0$. On the other hand, $\left(\frac{d}{dx_1}\right)_q^{i_1} \cdots \left(\frac{d}{dx_n}\right)_q^{i_n} x_1^{i_1} \cdots x_n^{i_n} = [i_1]_q \cdots [i_n]_q$. □

**Remark 5.5.** Let $X$ be $\mathbb{P}^n$, with fan given by $v_i = e_i$ for $i = 1, \ldots, n$, and $v_{n+1} = -e_1 - \cdots - e_n$. The generalized Rogers-Szegő polynomial $RS_{kD_{n+1}}(x;q)$ is the classical multivariate Rogers-Szegő polynomial

$$(8) \qquad RS_{k,n} = \sum_{i_0 + i_1 + \cdots + i_n = k} \begin{bmatrix} k \\ i_0, i_1, \ldots, i_n \end{bmatrix}_q x_1^{i_1} \cdots x_n^{i_n},$$

which is a specialization of a single-row Macdonald polynomial [Sz26, H97]. In the algebra generated by the $x_i$ and $\left(\frac{d}{dx_i}\right)_q$, there are operators

$$R_i := \sum_{l=0}^{n} e_{l+1}(x)(q - 1)^l \left(\frac{d}{dx_i}\right)_q^l.$$

Setting $L_i := \left(\frac{d}{dx_i}\right)_q$, we have the identities

$$R_i(RS_{k-1,n}) = RS_{k,n}, \quad L_i(RS_{k,n}) = [k]_q RS_{k-1,n}, \quad \text{and} \quad [L_i, R_i](RS_{k,n}) = q^k RS_{k,n}.$$

## 6. Asymptotics and an analogue of DH measure

We now prove some results about measures derived from $\chi_\mathbb{T}(\mathcal{Q}(\Delta)_\infty, \mathcal{O}(D^0))$. We start outside of the radially symmetric case, but specialize to it later. In particular, we will end with a limit theorem describing the asymptotic behavior of $RS_{kD}(x;q)$ as $k$ goes to infinity.

We again identify $M$ with $\mathbb{Z}^n$. For a Laurent polynomial in $x_i$ with $q$-series coefficients, $f(x) = \sum_u a_u(q) x^u$, its Fourier transform is defined as follows:

**Definition 6.1.** Let $\delta_u(y)$ be the Dirac measure at $u \in \mathbb{R}^n$. Let the measure $FT(f(x))$ be

$$FT(f(x)) = \sum_u a_u(q) \delta_u(y).$$

Note that if $g(x)$ is the function $\sum_u a_u(q) e^{-iu \cdot x}$ for $x \in \mathbb{R}^n$, and $q$ is chosen to that the $a_u(q)$ are real numbers, then $FT(f(x))$ is indeed the Fourier transform of $g(x)$.

Let $\tau_c : \mathbb{R}^n \to \mathbb{R}^n$ be dilation by $c$. For a $T = (\mathbb{C}^*)^n$-equivariant line bundle on a projective $T$-variety $Y$, Duistermaat-Heckman (DH) measures can be defined algebraically, following [BP90], as

$$(9) \qquad DH(Y, L) := \lim_{k \to \infty} (\tau_k)_* \left( \frac{FT\left(\chi_T(Y, L^{\otimes k})\right)}{k^{\dim X}} \right).$$

If $Y$ is instead an ind-projective ind-variety such as $\mathcal{Q}(\Delta)_\infty$, the measure defined above no longer makes sense. The factor $k^{\dim X}$ is not well-defined, and $FT(\chi_T(Y, L^{\otimes k}))$ may no longer properly define a distribution.

We let $Y = \mathcal{Q}(\Delta)_\infty$ be the quasimap ind-variety of a toric variety $X = X(\Delta)$, which has an additional $\mathbb{C}^*$-action that can be used to grade $\chi_T(Y, L^{\otimes k})$. Then, we define the following measure for a nef divisor $D$ on $X$.

**Definition 6.2.** Let $D$ be a $T$-invariant nef divisor on a smooth complete toric variety $X$. The *probability measure associated to $D$* is

$$(10) \qquad \mu_D = \lim_{q \to 1} \frac{FT(\chi_\mathbb{T}(\mathcal{Q}(\Delta)_\infty, \mathcal{O}(D^0)))}{\chi_\mathbb{T}(\mathcal{Q}(\Delta)_\infty, \mathcal{O}(D^0))|_{x=1}}.$$

In the radially symmetric case, this becomes

$$\mu_D = \lim_{q \to 1} \frac{FT(RS_D(x;q))}{RS_D(1;q)},$$

by Proposition 5.1.

It is easy to check that $\mu_D$ is a probability measure. The following proposition describes the support of $\mu_D$. We use the notation $v_\Delta$ for the sum of all primitive ray generators:
$$v_\Delta = \sum_i v_i.$$

**Proposition 6.3.** *Let $D = \sum_i a_i D_i$, and let $F$ be the face of the corresponding polytope $P = \bigcap_i \{u \,|\, \langle u, v_i \rangle \geq -a_i\}$ where the function $u \mapsto \langle u, v_\Delta \rangle$ is maximized. Then the measure $\mu_D$ is expressed by the formula*

$$\mu_D := \frac{\sum_{F \cap M} \binom{\langle u, v_\Delta \rangle + |a|}{\langle u, v_1 \rangle + a_1, \ldots, \langle u, v_r \rangle + a_r} \delta_u}{\sum_{F \cap M} \binom{\langle u, v_\Delta \rangle + |a|}{\langle u, v_1 \rangle + a_1, \ldots, \langle u, v_r \rangle + a_r}}. \qquad (11)$$

*In particular, $\mu_D$ is supported on $F$.*

*Proof.* By Theorem 4.1, we can rewrite $\mu_D$ as
$$\mu_D = \lim_{q \to 1} \sum_{u \in P \cap M} \frac{g_{\mathcal{A}_{P,u}}}{\sum_{u' \in P \cap M} g_{\mathcal{A}_{P,u'}}} \delta_u.$$

Let $u$ be in $P \cap M$. From (2),
$$g_{\mathcal{A}_{P,u}} = \prod_{i=1}^{r} \frac{1}{(q;q)_{\langle u, v_i \rangle + a_i}}.$$

If $u_F$ is any element of $F$, we have $c_F := \langle u_F, v_\Delta \rangle + |a| \geq \langle u, v_\Delta \rangle + |a|$, with equality precisely when $u \in F$. Thus,
$$(q;q)_{c_F} \cdot g_{\mathcal{A}_{P,u}} = \left( \prod_{l=1+\langle u, v_\Delta \rangle + |a|}^{c_F} (1-q^l) \right) \cdot \binom{\langle u, v_\Delta \rangle + |a|}{\langle u, v_i \rangle + a_i, \ldots, \langle u, v_r \rangle + a_r}_q.$$

The first factor on the right-hand side is 1 if $u \in F$, and vanishes at $q = 1$ otherwise, so
$$\lim_{q \to 1} (q;q)_{c_F} \cdot g_{\mathcal{A}_{P,u}}(q) = \begin{cases} \binom{\langle u, v_\Delta \rangle + |a|}{\langle u, v_i \rangle + a_i, \ldots, \langle u, v_r \rangle + a_r} & u \text{ in } F \cap M, \\ 0 & \text{otherwise.} \end{cases}$$

So,
$$\mu_D = \lim_{q \to 1} \sum_{u \in P \cap M} \frac{g_{\mathcal{A}_{P,u}}}{\sum_{u' \in P \cap M} g_{\mathcal{A}_{P,u'}}} \delta_u = \lim_{q \to 1} \sum_{u \in P \cap M} \frac{(q;q)_{c_F}}{(q;q)_{c_F}} \cdot \frac{g_{\mathcal{A}_{P,u}}}{\sum_{u' \in P \cap M} g_{\mathcal{A}_{P,u'}}} \delta_u,$$

which simplifies to the proposition. $\square$

Note that in the radially symmetric case, $F = P$.

**Example 6.4.** Let $X = \mathbb{P}^1$, and $D = [0]$. Then $P = \bigcap_i \{u \mid \langle u, v_i \rangle \geq -a_i\}$ is the interval $[0, 1]$ in $\mathbb{R}$. Then

$$\mu_{kD} = \frac{1}{2^k} \sum_{l=0}^{k} \binom{k}{l} \delta_l.$$

It is tempting to use $\mu_{kD}$ as a replacement of $\frac{FT(\chi_T(\mathcal{Q}(\Delta)_\infty, \mathcal{O}(kD^0)))}{k^{\dim \mathcal{Q}(\Delta)_\infty}}$ in the definition of DH measure. This produces a measure whose support must be contained in $P$, as for the usual DH measure. However, in the limit $\lim_{k \to \infty} (\tau_k)_* \mu_{kD}$, the measure concentrates too much to retain useful information about $P$, as illustrated in the following example.

**Example 6.5.** Let $D$ be as in Example 6.4. Then, let $\chi_{(\tau_k)_* \mu_{kD}}(x) = \int_\mathbb{R} e^{iyx} (\tau_k)_* \mu_{kD}(y)$ be the characteristic function (see e.g. [D10]) of the measure $(\tau_k)_* \mu_{kD}$. Explicitly,

$$\chi_{(\tau_k)_* \mu_{kD}}(x) = \frac{1}{2^k} \sum_{l=0}^{k} \binom{k}{l} e^{ilx/k} = \left(\frac{1 + e^{ix/k}}{2}\right)^k = e^{ix/2} \left(\frac{e^{-ix/2k} + e^{ix/2k}}{2}\right)^k = e^{ix/2} \cos^k(x/2k).$$

So $\lim_{k \to \infty} \chi_{(\tau_k)_* \mu_{kD}}(x) = e^{ix/2}$, and correspondingly $\lim_{k \to \infty} (\tau_k)_* \mu_{kD}(y) = \delta_{1/2}(y)$.

This motivates the following definition. Once again, $X = X(\Delta)$ is a smooth projective toric variety, $D$ a $T$-invariant nef divisor on $X$, and $D^0$ the corresponding divisor on the quasimap ind-variety $\mathcal{Q}(\Delta)_\infty$. We assume now that the generators $v_i$ of the rays of $\Delta$ sum to 0; in other words, $X$ has an ample divisor whose corresponding polytope is radially symmetric. For a measure $\mu$ on $M_\mathbb{R} \cong \mathbb{R}^n$, let $E_\mu = (\int_{\mathbb{R}^n} y_i d\mu)$ denote the corresponding expected value of the vector $(y_1, \ldots, y_n)$.

**Definition 6.6.** The *centered loop-DH measure* associated to $D$ is

$$\nu_{D,\text{cent}} = \lim_{k \to \infty} (\tau_{\sqrt{k}})_* \left(\mu_{kD} * \delta_{-E_{\mu_{kD}}}\right).$$

We define the *potential of D* as

$$\varphi_D(m) = \prod_{i=1}^{r} (\langle m, v_i \rangle + a_i)^{(\langle m, v_i \rangle + a_i)}.$$

On the Newton polytope of $D$, $\varphi_D$ has a unique minimum, which we call $m_D$. We define the *loop-DH measure* to be

$$\nu_D = \nu_{D,\text{cent}} * \delta_{m_D}.$$

The measures $\nu_{D,\text{cent}}$ and $\nu_D$ are Gaussian and can be expressed naturally in terms of the fan of $X$. Let $L_D$ be the linear subspace of $M_{\mathbb{R}}$ generated by differences of vectors in the Newton polytope of $D$.

**Theorem 6.7.** *Let $X$ be a smooth complete toric variety with a radially symmetric fan, and $D = \sum_i a_i D_i$ a $T$-invariant nef divisor. The probability measure $\nu_{D,\text{cent}}(u)$ on $M_{\mathbb{R}}$ is given by a Gaussian density function with mean $0$ and covariance matrix associated to the quadratic form*

$$u \mapsto \sum_{i \in I_D} \frac{1}{\langle m_D, v_i \rangle + a_i} \langle u, v_i \rangle^2,$$

*times the Dirac measure of $L_D$; that is, $\nu_{D,\text{cent}}(u)$ is a scalar multiple of*

$$(12) \qquad \exp\left(-\frac{1}{2} \sum_{i \in I_D} \frac{\langle u, v_i \rangle^2}{\langle m_D, v_i \rangle + a_i}\right) \cdot \delta_{L_D}(u).$$

*Here $I_D$ is the set of indices $i$ such that $\langle u, v_i \rangle + a_i$ is not uniformly $0$ on all of $P$.*

The proof is purely analytic and is given in the companion note [S23+]. Theorem 2 from the introduction follows, since convolving with $\delta_{m_D}$ replaces $u$ with $u - m_D$, and when $D$ is ample, $L_D = M_{\mathbb{R}}$.

To conclude, we consider a special case in which the proof that $\nu_D$ is Gaussian is very simple. It is inspired by results showing the pushforward in quantum $K$-theory of a homogeneous space becomes a ring homomorphism at $q = 1$ [BC18], although one should be warned that the $q$ considered by those authors plays a rather different role from the $q$ appearing here.

Let $\mathscr{M}$ denote the set of probability measures on $M_{\mathbb{R}}$, which is a semigroup with respect to convolution.

**Theorem 6.8.** *Suppose that $\mathcal{N} \subset \text{Div}_T(X)$ is a semigroup of divisors such that for $D, D' \in \mathcal{N}$,*

$$\frac{\chi_{\mathbb{T}}(\mathcal{Q}(\Delta)_\infty, \mathcal{O}(D^0 + D'^0))}{\chi_{\mathbb{T}}(\mathcal{Q}(\Delta)_\infty, \mathcal{O}(D^0 + D'^0))|_{x=1}}\bigg|_{q=1} = \frac{\chi_{\mathbb{T}}(\mathcal{Q}(\Delta)_\infty, \mathcal{O}(D^0))}{\chi_{\mathbb{T}}(\mathcal{Q}(\Delta)_\infty, \mathcal{O}(D^0))|_{x=1}}\bigg|_{q=1} \cdot \frac{\chi_{\mathbb{T}}(\mathcal{Q}(\Delta)_\infty, \mathcal{O}(D'^0))}{\chi_{\mathbb{T}}(\mathcal{Q}(\Delta)_\infty, \mathcal{O}(D'^0))|_{x=1}}\bigg|_{q=1},$$

*or equivalently, that $\mu_{(-)} : \mathcal{N} \to \mathscr{M}$ is a homomorphism. Then for each $D \in \mathcal{N}$, $\nu_D$ is also Gaussian, and*

$$\nu_{(-)} : \mathcal{N} \to \mathscr{M},$$

*is also a homomorphism, landing in the subset of Gaussian measures.*

*Proof.* We first verify that $\nu_{D,\text{cent}}$ is Gaussian. By definition,

$$\nu_{D,\text{cent}} = \lim_{k\to\infty}(\tau_{\sqrt{k}})_*\left(\mu_{kD} * \delta_{-E_{\mu_{kD}}}\right).$$

Since we have assumed that $\mu_{(-)}$ is a homomorphism, $\mu_{kD} = (\mu_D)^{*k}$. If $V_D$ is the random vector in $M_\mathbb{R}$ represented by the measure $\mu_D$, then the equation above means that $\mu_{kD}$ is the measure corresponding to a sum $\sum_{i=1}^k V_i$ of independent random vectors $V_i$, each with the same distribution as $V_D$. Thus $E_{\mu_{kD}} = kE(V_D)$ and the measure $\mu_{D,\text{cent}}$ corresponds to the random variable

$$\lim_{k\to\infty}\frac{\sum_{i=1}^k(V_i - E(V_D))}{\sqrt{k}},$$

which is a Gaussian with covariance equal to $Cov(V_D)$, by the central limit theorem as stated in [D10, Section 3.10].

We must also verify that $\nu_{D+D'} = \nu_D * \nu_{D'}$. First we check that $m_{D+D'} = m_D + m_{D'}$ under our assumptions on $D$ and $D'$. (It is not true in general!) In [S23+], it is shown that $\frac{E_{\mu_{kD}}}{k} \to m_D$. So $\frac{E_{\mu_{k(D+D')}}}{k} \to m_{D+D'}$. By assumption $\mu_{D+D'} = \mu_D * \mu_{D'}$, so we also have that $\frac{E_{\mu_{k(D+D')}}}{k} = \frac{E_{\mu_{kD}*\mu_{kD'}}}{k} = \frac{E_{\mu_{kD}}+E_{\mu_{kD'}}}{k} \to m_D + m_{D'}$. Thus $m_{D+D'} = m_D + m_{D'}$.

Finally, we calculate

$$\nu_{D+D'} = \lim_{k\to\infty}(\tau_{\sqrt{k}})_*\left(\mu_{k(D+D')} * \delta_{-E_{\mu_{k(D+D')}}}\right) * \delta_{m_{D+D'}}$$

$$= \lim_{k\to\infty}(\tau_{\sqrt{k}})_*\left(\mu_{kD} * \delta_{-kE_{\mu_D}} * \mu_{kD'} * \delta_{-kE_{\mu_{D'}}}\right) * \delta_{m_D} * \delta_{m_{D'}}$$

$$= \nu_D * \nu_{D'}$$

as desired, using $m_{D+D'} = m_D + m_{D'}$ in the second line. □

**Remark 6.9.** The above theorem applies for instance when $X = \mathbb{P}^n$. The primitive elements in the rays of the fan of $X$ are $v_i = e_i$ for $i = 1,\ldots,n$ and $v_{n+1} = -\sum_i e_i$. Then

$$\mu_{kD_{n+1}} = \frac{1}{(n+1)^k}\sum_{i_0+\cdots+i_n=k}\binom{k}{i_0,\ldots,i_n}\delta_{(i_1,\ldots,i_n)},$$

so $\mu_{kD_{n+1}} * \mu_{lD_{n+1}} = \mu_{(k+l)D_{n+1}}$.

<div style="text-align:center">REFERENCES</div>

*Email address*: `anderson.2804@math.osu.edu`

Department of Mathematics, The Ohio State University, Columbus, OH 43210, USA

*Email address*: `shah@karlin.mff.cuni.cz`

Department of Algebra, Faculty of Mathematics and Physics, Charles University, Prague, Czech Republic